\documentclass[10pt]{article}
\usepackage[utf8]{inputenc}
\usepackage[T1]{fontenc}
\usepackage{amsmath}
\usepackage{amsfonts}
\usepackage{amssymb}
\usepackage{stmaryrd}
\usepackage{amsthm}
\usepackage{graphicx}
\usepackage{hyperref}
\usepackage{enumitem}

\DeclareMathOperator{\diag}{diag}
\newcommand{\subscript}[2]{$#1 _ #2$}

\graphicspath{figures/}

\usepackage{authblk}
\usepackage[
backend=biber,
sorting=none
]{biblatex}
\title{ Optimal Real Time Drone Path  Planning for Harvesting Information from a Wireless Sensor Network}
\author[1]{Ramkumar Ganapathy}
\author[1]{Christopher Thron\thanks{Corresponding author: thron@tamuct.edu}}
\affil[1]{Texas A\&M University-Central Texas, Killeen TX USA}

\date{}

\begin{document}
\maketitle

\begin{abstract}

We consider a remote sensing system in which fixed sensors are placed in a region, and a drone flies over the region to collect information from cluster heads. We assume that the drone has a fixed maximum range,  and that the energy consumption for information transmission from the cluster heads increases with distance according to a power law. Given these assumptions, we derive local optimum conditions for a drone path that either minimizes the total energy  or the maximum energy required by the cluster heads to transmit information to the drone.  We show how a homotopy approach can produce a family of solutions for different drone path lengths, so that a locally optimal solution can be found for any drone range. We implement the homotopy solution in python, and demonstrate the tradeoff between drone range and cluster head power consumption for several geometries. Execution time is sufficiently rapid for the computation to be performed real time, so the drone path can be recalculated on the fly. The solution is shown to be globally optimal for sufficiently long drone path lengths. For future work, we indicate how the solution can be modified to accommodate moving sensors.
\end{abstract}

Keywords: Wireless sensor network, drone, path planning, information collection, power-efficient, extended lifetime, optimization, 


\section{Introduction}\label{sec:intro}
There are many critical applications  for remote sensing in low-resourced areas. Some of these include agricultural monitoring of crops and/or forests for fire and/or disease; wildlife monitoring, to track wildlife movement and detect poachers; free-range livestock monitoring, to track herd movement and to prevent cattle rustling; and early warning of terrorist or bandit activity.  In low-resourced situations especially, the system cost is of huge importance, and spells the difference as to whether the system may be implemented or not. Sensors in area-monitoring networks typically are battery-powered, and must be replaced when their batteries is exhausted. This can be both difficult and costly, particularly in inaccessible locations. For this reason, reducing the power consumption of sensors in the field is a key factor in designing such remote sensing systems.

Previous authors have investigated practical use cases for using drones to collect information from WSN's. \cite{nguyen2021uav} has provided a comprehensive survey of drone-assisted data collection from wireless sensor networks. Several previous papers have dealt specifically with drone path planning in such scenarios. \cite{ho2013heuristic} gives a heuristic algorithm to  decide which node within a cluster to fly to in order to gather information from the cluster. \cite{liu2018performance} uses a boustrophedon-type  flight path for a sensor network located in region divided into square cells. \cite{gong2018flight} considers the problem of minimizing flight time of an information-gathering drone in the case where sensors are in a straight line, and the time required for information transfer from each sensor is an important consideration. \cite{skiadopoulos2020impact} considers the impact of drone path shape on information transmission from sensors to drone.
\cite{6825191} uses particle-swarm optimization to arrive at an optimal path.

\begin{figure}
    \centering
    \includegraphics[width=4.5 in]{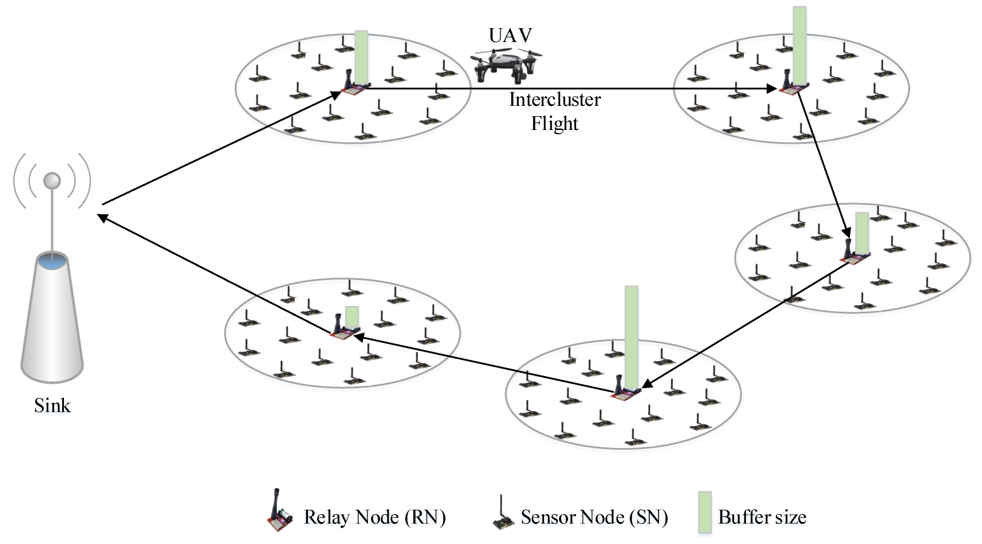}
    \caption{Schematic of remote sensing system in which a drone harvests information from fixed sensors (copied from Reference \cite{nguyen2021uav}).}
    \label{fig:schematic}
\end{figure}

In this paper we consider a hybrid WSN-drone system for remote sensing. A single drone is used to collect information from wireless sensors which have fixed locations in the field (see Figure~\ref{fig:schematic}.  The drone harvests information from all sensors during a single tour.  The system design is posed as a constrained optimization problem. Our novel solution approach involves using a Lagrange multiplier to construct a differential equation in which the Lagrange multiplier is the independent variable. Our approach has reduced computational complexity, so that real-time solution within seconds is possible even for very large systems. 

\section{Methodology}\label{sec:method}
\subsection{System assumptions}
We consider a remote sensing system that satisfies the following assumptions:
\begin{enumerate}[label=(\subscript{A}{{\arabic*}})]
\item
The system includes a wireless sensor network with fixed cluster heads.  Each sensor in the field transmits its information (either directly or indirectly) to one of the cluster heads. 
\item
The system also includes a drone that flies on a specific trajectory (to be determined), 
so that each cluster head transmits its energy to the drone when it is nearby. 
\item
The drone's energy consumption is determined entirely by the length of the drone's trajectory. (This is the case for example if the drone flies at constant velocity, and there is no wind or other conditions that would affect the flight of the drone.)
\item
The drone has a fixed maximum path length, which is determined by the energy capacity of the battery and the rate of energy consumption. 
\item
The drone's  location is approximately constant while receiving information from a particular cluster head.  This assumption is satisfied in either of the following scenarios: (a) the information transmission from the cluster head occupies only a brief period of time, so that the distance the drone moves during the transmission period is negligible; or (b) the drone lands during transmission, and consumes no energy during the transmission period.  
\item 
The energy expended by the cluster head in information transmission is proportional to the distance between the cluster head and drone raised to a power law exponent.
\end{enumerate}
We also consider two alternative criteria for optimization:
\begin{enumerate}[label=(\subscript{O}{{\arabic*}})]
\item Minimize the total energy expended by cluster heads in information transmission;
\item Minimize the maximum of the energies expended by cluster heads in information transmission.
\end{enumerate}

Under assumptions $A_1$-$A_6$ and criteria $O_1$-$O_2$, it follows that an optimum path for the drone will  consist of straight-line segments joining the drones' locations where it harvests energy from the different cluster heads. 
This conclusion is reflected in the mathematical description in subsequent sections.

\subsection{System parameters and variables}
The system parameters include:
\begin{itemize}
\item
$(x_{j},y_{j}), j = 1...J$: positions of the cluster heads that are sending the information;
\item 
$L$: Maximum path length for the drone;
\item 
$p$: Power loss exponent.
\end{itemize}

\noindent
The system variables include:
\begin{itemize}
\item 
$(u_{j},v_{j}), j = 1...J$: Drone positions for harvesting information from cluster heads. Here the order of points the drone takes is assumed to be $(u_{j}, v_{j})$ to $(u_{j+1}, v_{j+1})$ for $j = 1\ldots J-1$. 
\item
$(u_0,v_0)$, $(u_{J+1},v_{J+1})$: Drone starting and ending position, respectively.
\end{itemize}

\subsection{Mathematical formulation of the optimization problem}
Information transmission from cluster head $(x_j,y_j)$  takes place when the drone is at position $(u_j,v_j)$. The total energy consumption required for information transfer is the sum of the energies from the $J$ cluster heads, which following assumption $(A_6)$ can be  represented as a function $f(\vec{u}, \vec{v})$ where:

\begin{equation}\label{eq:def_f}
f(\vec{u}, \vec{v}) = \sum_{j=1}^J \left((x_j-u_j)^2 + (y_j - v_j)^2\right)^{p/2},
\end{equation}
where $\vec{u}:= [u_1,\ldots u_J]$ and $\vec{v}:= [v_1,\ldots v_J]$.

The distance the drone travels is given by $g(\vec{u}, \vec{v})$, where
\begin{equation}\label{eq:def_g}
g(\vec{u}, \vec{v}) = \sum_{j=0}^J \sqrt{ (u_j-u_{j+1})^2 + (v_j-v_{j+1})^2 }
\end{equation}

The system's constraint is given by $g(\vec{u}, \vec{v})\le L$, where $L$ is the maximum distance the drone can travel. In the case where $L$ is large enough so that the drone is capable of making a complete tour of the cluster heads, then any such complete tour will give $g(\vec{u},\vec{v}) = 0$ and is thus an optimal solution to the problem.  So we consider instead the more difficult problem where $L$ is smaller than the smallest tour distance, so $g(\vec{u}, \vec{v}) = 0$ is impossible. 
In this case we may replace the inequality with an equality constraint. This can be seen as follows. Suppose a drone path includes the point $(u_j,v_j)$ where $|(u_j,v_j) - (x_j,y_j)|>0$. For simplicity we define: 
\begin{equation}\label{eq:wzDef}
\vec{w}_j := (u_j,v_j) \text{~and~}  \vec{z}_j := (x_j,y_j).    
\end{equation}
 Then for any $\delta$ with $0<\delta \le 1$ we have 
\begin{equation}
|\vec{w}_j + \delta(\vec{z}_j - \vec{w}_j) - \vec{z}_j| = (1-\delta)|\vec{w}_j - \vec{z}_j| < |\vec{w}_j - \vec{z}_j|,
\end{equation}
so replacing $\vec{w}_j$ with $\vec{w}_j + \delta(\vec{z}_j - \vec{w}_j)$ will reduce the energy function \eqref{eq:def_f}.  
Now if the drone path length is less than $L$, then $\delta>0$ can be chosen sufficiently small such that replacing  $\vec{w}_j$ with $\vec{w}_j + \delta(\vec{z}_j - \vec{w}_j)$ still yields total drone path length less than $L$. Thus any drone path with length less than $L$ can be improved: so no drone path with length less than $L$ can be optimal. 

We may now formulate the optimization problem. In case $O_1$, we have
\begin{equation}\label{eq:fundOpt}
\begin{aligned}
\text{Minimize } f(\vec{u},\vec{v})
\text{  subject to }
g(\vec{u},\vec{v}) = L.
\end{aligned}
\end{equation}
In case $O_2$, minimizing the maximum transmission energy is equivalent to minimizing $|\vec{w}_j -\vec{z}_J|$ over all $j$.  Using the equality
\begin{equation}
\lim_{p\rightarrow \infty} (a_1^p +  \ldots +  a_J^p)^{1/p} = \max_j a_j\quad \text{ if }a_j>0~\forall j, 
\end{equation}
with $a_j := |\vec{w}_j -\vec{z}_j|$, we conclude that by taking $p$ sufficiently large we can approach the solution to $O_2$ with arbitrary precision.
\subsection{Local optimization conditions}\label{sec:locOpt}
The minimization problem \eqref{eq:fundOpt} leads to the following Lagrange multiplier condition for a local minimum:
\begin{equation}\label{eq:Lagrange}
\nabla f(\vec{u}, \vec{v}) = \lambda \nabla(g(\vec{u}, \vec{v}) - L),
\end{equation}
where $\nabla h$ denotes the gradient of $h$ with respect to $\vec{u},\vec{v}$:
\begin{equation}
\nabla h(\vec{u}, \vec{v}) := \left(\frac{\partial f}{\partial u_1},  \ldots, \frac{\partial f}{\partial u_J}, \frac{\partial f}{\partial v_1},\ldots,\frac{\partial f}{\partial v_J}\right),
\end{equation}
and $\lambda$ is a Lagrange multiplier. 

Writing the vector equation \eqref{eq:Lagrange} out in terms of components, we have:

\begin{equation}\label{eq:opt1}
\frac{\partial f}{\partial u_j} = \lambda\frac{\partial g}{\partial u_j};\qquad
\frac{\partial f}{\partial v_j} = \lambda\frac{\partial g}{\partial v_j} 
\qquad (j=1 \ldots J).
\end{equation}

For notational simplicity we define new variables. Let $(j=1\ldots J)$
\begin{equation}\label{eq:def_simp_nots}
\begin{aligned}
&{a_j} := {u_j - x_j}; \qquad &{b_j} := {v_j - y_j}; \qquad 
&A_j := {a_j}^2 + {b_j}^2, 
\end{aligned}
\end{equation}
and let $(j=0\ldots J)$
\begin{equation}\label{eq:def_simp_nots2}
\begin{aligned}
&{m_j} = u_{j} - u_{j+1}; \qquad &{n_j} := v_{j} - v_{j+1};\qquad 
&M_j := ({m_j}^2 + {n_j}^2)^{1/2}.
\end{aligned}
\end{equation}
Then we have:
\begin{equation}\label{eq:def_f_g}
f(\vec{u}, \vec{v}) = \sum_{j=1}^J A_j^{p/2}~~\text{and}~~
g(\vec{u}, \vec{v}) = \sum_{j=0}^J M_j.
\end{equation}
Using this notation, we have:
\begin{equation}
\begin{aligned}
\frac{\partial f}{\partial u_j} &= \frac{\partial f}{\partial A_j} \cdot\frac{\partial A_j}{\partial a_j} \cdot \frac{\partial a_j}{\partial u_j}
= ((p/2)A_j^{p/2-1} \cdot 2a_j \\
&= pa_j A_j^q,~~ \text{where }q:=p/2-1,)
\end{aligned}
\end{equation}
and
\begin{equation}\label{eq:pdu}
\begin{aligned}
\frac{\partial g}{\partial u_j} &= \frac{\partial g}{\partial M_{j-1}} \cdot \frac{\partial M_{j-1}}{\partial m_{j-1}} \cdot \frac{\partial m_{j-1}}{\partial u_j} + \frac{\partial g}{\partial M_j} \cdot \frac{\partial M_j}{\partial m_j} \cdot \frac{\partial m_j}{\partial u_j}\\
&= (1/2)M_{j-1}^{-1}(2m_{j-1})(-1) + (1/2)M_j^{-1}(2m_j)\\
&= \frac{m_{j}}{M_{j}}  - \frac{m_{j-1}}{M_{j-1}}.
\end{aligned}
\end{equation}
Similarly we have: 
\begin{equation}\label{eq:pdv}
\frac{\partial f}{\partial v_j} = pb_j A_j^q;\qquad 
\frac{\partial g}{\partial v_j} = \frac{n_{j}}{M_{j}}  - \frac{n_{j-1}}{M_{j-1}}.
\end{equation}
It follows that we may rewrite the equations in \eqref{eq:opt1} as:
\begin{equation}\label{eq:opt2}
pa_j A_j^q = \lambda\left( \frac{m_{j}}{M_{j}}  - \frac{m_{j-1}}{M_{j-1}}\right); \quad
pb_j A_j^q = \lambda\left( \frac{n_{j}}{M_{j}}  - \frac{n_{j-1}}{M_{j-1}}\right)~~(j=1\ldots J).
\end{equation}
Equations \eqref{eq:opt2} are required for local optimality. Hence they are necessary for global optimality, but not sufficient. So the solutions that we provide in this paper may or may not be globally optimal. We return to this issue in Section~\ref{sec:future}. 
\subsection{Homotopy approach to local optimization for a given path length}\label{sec:homotopy}

Solutions of the fundamental optimization problem \eqref{eq:fundOpt}
 for different values of $L$ correspond to solutions of \eqref{eq:opt2} with different values of $\lambda$.  As $L$ is varied continuously, the value of $\lambda$ will also vary continuously, and the components of $\vec{u}$ and $\vec{v}$ will vary continuously as well. This suggests that if we have a solution to  \eqref{eq:fundOpt} for a given value of $L$, then we may be able to perturb that solution differentially to obtain solutions for different values of $L$. In fact, we do have a solution for a particular value of $L$: namely $f(\vec{u},\vec{v})=0$ when $u_j=x_j, v_j=y_j$, which solves the constraint  $g(\vec{u},\vec{v}) = L$ when $L$ is  equal to the length of the  shortest tour that visits all sensors (the so-called ``travelling salesman'' tour). 
 
 In practice, we want to find solutions corresponding to  values of $L$ that are smaller than the length of a full tour. So we start with the full tour, and successively ``nudge'' the solutions so they correspond to smaller and smaller values of $L$. This is an example of a \emph{homotopy approach}: start with a solution that is optimal for a different set of conditions, and generate a smooth curve of solutions that joins this solution to a solution that satisfies desired conditions. In practice this smooth curve of solutions is often specified parametrically as the solution to  a set of ordinary differential equations in a given parameter.  In the next section, we will derive differential equations that can be used to find the parametrized curve that leads to our desired optimal solution.

\subsection{Derivation of differential equations for parametrized homotopy curve}\label{sec:derivation}
 
 First  we parametrize solutions to \eqref{eq:opt2} using the parameter $s$, so that $a_j, b_j, A_j, m_j, M_j,$ \text{ and } $\lambda$ are all functions of $s$.  As $s$ changes the equalities cannot change, so we have
\begin{equation}\label{eq:optDiff}
\begin{aligned}
p \frac{d}{ds} \left(a_j A_j^q\right) = \frac{d}{ds} \left( \lambda\frac{m_{j}}{M_{j}}  - \lambda\frac{m_{j-1}}{M_{j-1}}\right) \\
p \frac{d}{ds} \left(b_j A_j^q\right) = \frac{d}{ds} \left( \lambda\frac{n_{j}}{M_{j}}  - \lambda\frac{n_{j-1}}{M_{j-1}}\right)
\end{aligned}
\end{equation}

 Since we want to find solutions for smaller values of $L$, we want to ensure that $L$ decreases as $s$ increases.  We may ensure this by adding the condition:
 \begin{equation}\label{eq:gDeriv}
 \frac{dg}{ds} = -1 \implies \sum_{j=0}^J \frac{dM_j}{ds} = -1.
 \end{equation}
Recall that the variables $a_j, b_j, A_j, m_j, M_j,$ in \eqref{eq:optDiff}
 are all functions of $\vec{u}$ and $\vec{v}$, so the derivatives in \eqref{eq:optDiff} can all be expressed in terms of the derivatives $\frac{du_j}{ds}$ and $\frac{dv_j}{ds}$ for $j=1,\ldots J$. As a result, we obtain a system of $2J+1$ differential equations for the variables $\{u_j,v_j\}, j=1\ldots J$ and $\lambda$. 

We may rewrite the first equation in \eqref{eq:optDiff} as $(j=1\ldots J)$:
\begin{equation}\label{eq:expandEq}
\begin{aligned}
0 = &p  A_j^q \frac{da_j}{ds} + pqa_jA_j^{q-1}\frac{dA_j}{ds} -\frac{\lambda} {M_j}\frac{dm_j}{ds} + \frac{\lambda}{ M_{j-1}}\frac{dm_{j-1}}{ds}\\
&+\frac{\lambda m_j} {M_j^2}\frac{dM_j}{ds} -\frac{\lambda m_{j-1}}{M_{j-1}^2}\frac{dM_{j-1}}{ds}
-\frac{m_j}{M_j}\frac{d\lambda}{ds} +\frac{m_{j-1}}{M_{j-1}}\frac{d\lambda}{ds}
\end{aligned}
\end{equation}
 The derivatives in \eqref{eq:expandEq} may be expressed in terms of $\left\{\frac{du_k}{ds}\right\}$ and $\left\{\frac{dv_j}{ds}\right\}$ using the formulas $(j=0\ldots J)$::
 \begin{equation}
     \begin{aligned}
         \frac{da_j}{ds} &= \frac{du_j}{ds};\\
         \frac{dA_j}{ds} &= 2a_j\frac{du_j}{ds} + 2b_j\frac{dv_j}{ds};\\
         \frac{dm_j}{ds} &= \frac{du_j}{ds} - \frac{du_{j+1}}{ds};\\
         \frac{dM_j}{ds} &= M_j^{-1} \left( m_j \left(\frac{du_j}{ds} - \frac{du_{j+1}}{ds}\right) + n_j \left(\frac{dv_j}{ds} - \frac{dv_{j+1}}{ds}\right)\right);\\
     \end{aligned}
 \end{equation}
 It follows that \eqref{eq:expandEq} with index $j$ will depend on the derivatives $\frac{du_k}{ds}$, $\frac{dv_k}{ds}$ and $\frac{d\lambda}{ds}$ for $k=j-1,j,j+1$ (noting that $\frac{du_0}{ds} = \frac{dv_0}{ds} = \frac{du_{J+1}}{ds} = \frac{dv_{J+1}}{ds} = 0$). 
 By collecting terms  corresponding to each derivative, we find $(j=1\ldots J)$::
 \begin{equation}\label{eq:expandEq1}
\begin{aligned}
0 = &\frac{du_j}{ds}\left(p  A_j^q  + pqa_jA_j^{q-1}(2a_j) -\frac{\lambda} {M_j} - \frac{\lambda}{ M_{j-1}}
+\frac{\lambda m_j^2} {M_j^3} +\frac{\lambda m_{j-1}^2}{M_{j-1}^3}\right)\\
&+\frac{du_{j-1}}{ds}\left( \frac{\lambda} {M_{j-1}} -\frac{\lambda m_{j-1}^2} {M_{j-1}^3}\right)+\frac{du_{j+1}}{ds}\left( \frac{\lambda} {M_j} -\frac{\lambda m_j^2} {M_j^3}\right)\\
&+\frac{dv_j}{ds}\left(pqa_jA_j^{q-1}(2b_j)
+\frac{\lambda m_jn_j} {M_j^3} +\frac{\lambda m_{j-1}n_{j-1}}{M_{j-1}^3}\right)\\
&+\frac{dv_{j-1}}{ds}\left(-\frac{\lambda m_{j-1}n_{j-1}} {M_{j-1}^3}\right)+\frac{dv_{j+1}}{ds}\left(  -\frac{\lambda m_j n_j} {M_j^3}\right)\\
&+ \frac{d\lambda}{ds}\left(-\frac{ m_{j}} {M_{j}} + \frac{m_{j-1}} {M_{j-1}}\right)\\
\end{aligned}
\end{equation}
Using various algebraic manipulations and the fact that:
\begin{equation}
    1-\frac{m_j^2}{M_j^2} = \frac{n_j^2}{M_j^2},
\end{equation}
 the equations  \eqref{eq:expandEq} may be simplified to:
\begin{equation}\label{eq:expandEq2}
\begin{aligned}
0 = &\frac{du_j}{ds}\left(p  A_j^q\left(1 + \frac{2qa_j^2}{A_j}\right) - \frac{\lambda n_{j-1}^2}{ M_{j-1}^3} -\frac{\lambda n_j^2} {M_j^3} \right)\\
&+\frac{du_{j-1}}{ds}\left( \frac{\lambda n_{j-1}^2} {M_{j-1}^3} \right)+\frac{du_{j+1}}{ds}\left( \frac{\lambda n_j^2} {M_j^3} \right)\\
&+\frac{dv_j}{ds}\left(2pqa_jb_jA_j^{q-1}
+ \frac{\lambda m_{j-1}n_{j-1}}{M_{j-1}^3}
+\frac{\lambda m_jn_j} {M_j^3} \right)\\
&-\frac{dv_{j-1}}{ds}\left(\frac{\lambda m_{j-1}n_{j-1}} {M_{j-1}^3}\right)-\frac{dv_{j+1}}{ds}\left(  \frac{\lambda m_j n_j} {M_j^3}\right)\\
&+ \frac{d\lambda}{ds}\left(-\frac{ m_{j}} {M_{j}} + \frac{m_{j-1}} {M_{j-1}}\right), \quad j=1\ldots J:
\end{aligned}
\end{equation}
The  $J$ equations associated with the second equation in \eqref{eq:optDiff} may be obtained from \eqref{eq:expandEq2} by exchanging $u_k\leftrightarrow v_k$, $a_k\leftrightarrow b_k$, and $m_k\leftrightarrow n_k$. 

Finally, \eqref{eq:gDeriv}
can be expressed in terms of derivatives of $\{u_j,v_j\}$ as:
\begin{equation}
\sum_j \left(\frac{m_j}{M_j} - \frac{m_{j-1}}{M_{j-1}}\right)\frac{du_j}{ds}  + \left(\frac{n_j}{M_j} - \frac{n_{j-1}}{M_{j-1}}\right)\frac{dv_j}{ds}=-1.
\end{equation}
 
The entire system  of $2J+1$ equations  can be represented in matrix form as: 
\begin{equation}\label{eq:H}
H\frac{d\mathbf{u}}{ds} = \mathbf{z} ~\implies \frac{d\mathbf{u}}{ds} = H^{-1}\mathbf{z},
\end{equation}
Where $H$ is a $(2J+1)\times(2J+1)$ matrix,
and $\mathbf{u},\mathbf{z}$ are column vectors of length $2J+1$, given by the formulas:
 \begin{equation}
 \begin{aligned}
 \mathbf{u} &:= \left[\vec{u}, \vec{v}, \lambda\right]^T;\\
 \mathbf{z} &:= \left[0,\ldots,0,-1\right]^T
 \end{aligned}
 \end{equation}
(i.e. $\mathbf{z}$ has a single nonzero entry of -1 in the last component).

The matrix $H$ can be decomposed into blocks:
\begin{equation}\label{eq:Hblock}
H =  
\begin{bmatrix}
H_{11} & H_{12} & \vec{h_1}\\
H_{21} & H_{22} & \vec{h_2}\\
\vec{h}_1^T & \vec{h}_2^T & 0
\end{bmatrix},
\end{equation}
The blocks $H_{ij}$ can be expressed in terms of diagonal and subdiagonal matrices. First we introduce some abbreviated notations. Let $\diag(c_j)$ denote the $J \times J$ diagonal matrix with entries $c_1,\ldots c_J$ on the diagonal. Let $D$ denote the discrete forward derivative matrix with $-1$'s on the diagonal and 1's on the first superdiagonal. Then we have:
\begin{equation}\label{eq:Hspec}
\begin{aligned}
    H_{11} &= \diag\left( pA_j^q\left(1 + \frac{2qa_j^2}{A_j}\right) \right) + \lambda \diag \left( \frac{n_j^2}{M_j^3} \right) D + \lambda \diag \left( \frac{n_{j-1}^2}{M_{j-1}^3} \right)D^T\\
     H_{12} &= \diag\left( 2pqa_jb_jA_j^{q-1}\right) - \lambda \diag \left( \frac{m_jn_j}{M_j^3} \right) D - \lambda \diag \left( \frac{m_{j-1}n_{j-1}}{M_{j-1}}\right)D^T \\
     H_{21} &= H_{12}\\
     H_{22} &= \diag\left( pA_j^q\left(1 + \frac{2qb_j^2}{A_j}\right) \right) + \lambda \diag \left( \frac{m_j^2}{M_j^3} \right) D + \lambda \diag \left( \frac{m_{j-1}^2}{M_{j-1}^3} \right)D^T\\
     \vec{h}_1 &= \left(-\frac{m_1}{M_1} + \frac{m_{0}}{M_{0}}~,~ \ldots~,~ -\frac{m_J}{M_J} + \frac{m_{J-1}}{M_{J-1}} \right)^T \\
     \vec{h}_2 &= \left(-\frac{n_1}{M_1} + \frac{n_{0}}{M_{0}}~,~ \ldots ~,~-\frac{n_J}{M_J} + \frac{n_{J-1}}{M_{J-1}} \right)^T.
\end{aligned}
\end{equation}

\subsection{Initial conditions for differential equation}\label{sec:initialConditions}

In Section \ref{sec:homotopy}, we suggested that by starting with a complete tour by the drone that visits all sensors and ``nudging'' the drone's path, we can find optimal solutions for shorter and shorter path lengths.  Unfortunately, the matrix $H$ in \eqref{eq:H} is singular when $(u_j,v_j) = (x_j,y_j), j=1 \ldots J$, which corresponds exactly to the case of a complete tour.  So we cannot use this solution as an initial condition.  Fortunately we can address this problem by choosing a solution where the $(u_j,v_j)$  are a small offset from $ (x_j,y_j)$ for all $j$, i.e.
\begin{equation}\label{eq:uvInit}
    (u_j,v_j) = (x_j,y_j) + \vec{\epsilon}_j.
\end{equation}
We need to choose $\vec{\epsilon}_j$ in such a way that the Lagrange multiplier conditions \eqref{eq:Lagrange}  are satisfied. We also need to choose $\vec{\epsilon}_j$ such that the path joining the $ (u_j,v_j)$ is smaller than the full tour.



The gradient of $f(\vec{u},\vec{v})$ evaluated at the point for the initial point given by \eqref{eq:uvInit} is given by:

\begin{equation}\label{eq:initial_cond}
\begin{aligned}
\nabla{f(u_{j}, v_{j})} &= p\left(x_{j}-u_{j})^{2} + (y_{j}-v_{j})^{2}\right)^{(\frac{p}{2}-1)}\bigl[(x_{j}-u_{j}), (y_{j}-v_{j})\bigr]\\
&= p\left| \vec{\epsilon}_j \right|^{p-2}\vec{\epsilon}_j\\
&= p\left| \vec{\epsilon}_j \right|^{p-1}\hat{\epsilon}_j,
\end{aligned}
\end{equation}
where $\hat{\epsilon}_j$ is the unit vector in the direction of $\vec{\epsilon}$.  The Lagrange multiplier condition \eqref{eq:Lagrange}
 gives:
\begin{equation}\label{eq:eps_g}
\begin{aligned}
p\left| \vec{\epsilon}_j \right|^{p-1}\hat{\epsilon}_j = \lambda  \nabla_j g(\vec{u},\vec{v}),
\end{aligned}
\end{equation}
where 
\begin{equation}
    \nabla_j g(\vec{u},\vec{v}) := \left( \frac{\partial g}{\partial u_j}, \frac{\partial g}{\partial v_j} \right).
\end{equation}
Since $g(\vec{u},\vec{v})$ is smooth in the vicinity of $(\vec{u},\vec{v})\approx  (\vec{x},\vec{y})$, we may approximate
\begin{equation}\label{eq:gj}
\nabla_j g(\vec{u},\vec{v}) \approx 
\left. \nabla_j g(\vec{u},\vec{v}) \right|_{\vec{u} = \vec{x},\vec{v} = \vec{y}},
\end{equation}
which may be evaluated using \eqref{eq:pdu} and \eqref{eq:pdv}.  Denoting the right-hand side of \eqref{eq:gj} as $\vec{g}_j$, we have from \eqref{eq:eps_g}
\begin{equation}\label{eq:eps_g2}
\begin{aligned}
&p\left| \vec{\epsilon}_j \right|^{p-1}\hat{\epsilon}_j = \lambda  \vec{g}_j\\
\implies & \hat{\epsilon}_j = \pm \hat{g}_j \text{ and } \left|\epsilon_{j}\right| = \left|\frac{\lambda}{p}\right|^{(\frac{1}{p-1})}|\vec{g}_{j}|^{(\frac{1}{p-1})}.
\end{aligned}
\end{equation}
Since we are interested in solutions for which $g$ is reduced, we choose the negative sign in \eqref{eq:eps_g2}.  Then $\vec{\epsilon}_j$ is uniquely determined by the values of $\vec{g}_j$ and $\lambda$.  By choosing a small value of $\lambda$, we may obtain initial values for $(\vec{u},\vec{v})$ that are very close to $(\vec{x},\vec{y})$, so that the approximation in \eqref{eq:gj} holds to a very high degree of accuracy.  In this way, we obtain initial conditions for our differential equation \eqref{eq:H} for which $H$ is not singular, and can thus apply standard numerical methods for solution.

\subsection{Global optimality}

In Section~\ref{sec:locOpt} we posed the local optimality condition \eqref{eq:Lagrange} from which the vector differential equation \eqref{eq:H} is derived.  To show that a locally optimal solution is globally optimal, we would need to show that it improves over all other locally optimal solutions. We may establish this for the solutions described in Sections~\ref{sec:homotopy}-\ref{sec:initialConditions}.

Every local optimal solution for a given drone path length will be part of a homotopy of solutions that may be parametrized by path length. As path length increases in the homotopy, the energy consumption decreases until a minimum energy consumption is reached for the entire homotopy. This minimum energy is necessarily nonnegative--and if it is 0, then the homotopy must terminate at a tour that joins all of the cluster heads. The shortest cluster head tour (i.e. the traveling salesman solution) will be the unique optimum in the case when the drone path length is equal to this shortest tour.  It follows that for sufficiently long drone path lengths, the solution that belongs to the homotopy that includes the traveling salesman solution will be the unique global optimum.  This justifies our claim that the solution of \eqref{eq:fundOpt} calculated using the homotopy approach outlined in Sections~\ref{sec:homotopy}-\ref{sec:initialConditions} is optimal for sufficiently long drone path lengths. 

\subsection{Implementation in python}
\label{sec:implementation}

As described above the algorithm has two phases. First, an initial ordering of cluster heads is determined, such that the length of a tour joining the cluster heads in this order is minimized among all possible orderings;
Once the ordering is set, the homotopy of solutions is computed, using the initial conditions and system of differential equations described in Sections~\ref{sec:initialConditions} and \ref{sec:homotopy}.
The final path with desired path length is selected from the homotopy. 

For the first phase,  the shortest tour joining cluster heads is found using travelling salesman algorithm as implemented in the  \texttt{python-tsp} package.. This gives a set of ordered cluster head points, that is, ($x_{j}$,$y_{j}$) re-arranged appropriately. 

For the second phase, the homotopy solution is computed according to the following steps:\begin{itemize}
\item[Step 0:] Initialize  cluster head locations, initial value of $\lambda$, and step size $h$ (typically on the order of 0.1),
\item[Step 1:] Find a minimal initial path joining the  cluster head points  using travelling salesman algorithm.
\item[Step 2:] Determine the initial points of closest approach to be a small distance away from the cluster head points, according to the algorithm described in Section~\ref{sec:initialConditions}.
\item[Step 3:] Solve the homotopy equations \eqref{eq:H} numerically. Initially we used the \texttt{odeint} solver from the \texttt{scipy} package, but encountered  instabilities when the the path vertices began to merge. We had better success using a Runge Kutta-4 code obtained from the web and modified for our purpose.
\end{itemize}

The Appendix gives a more detailed description of the code as well as a link to a Github page where the code can be downloaded.

\subsection{Specification of test cases}\label{sec:testCases}
The code was tested with the following configurations:
\medskip 

\noindent Case 1:

Cluster head positions: [(2, 1), (2, 4), (6, 4), (6, 1)] 

Starting position of the drone: (0, 0) 

Transmission power loss exponent: 2
\medskip 

\noindent Case 2:

Cluster head positions: [(2, 1), (2, 4), (8, 2), (6, 4), (6, 1)] 

Starting position of the drone: (0, 0)

Transmission power loss exponent: 2
\medskip

\noindent
Case 3:

Cluster head positions: [(2, 1), (2, 4), (8, 2), (6, 4), (6, 1)]

Starting position of the drone: (3, 1) 
Transmission power loss exponent: 2
\medskip

\noindent
Case 4:

Cluster head positions: [(2, 1), (2, 4), (8, 2), (6, 4), (6, 1), (7, 3.5), (1, 2.5)]

Starting position of the drone: (0, 0)

Transmission power loss exponent: 2
\medskip

All cases used a step size of 0.1. The execution time was extremely rapid, and all scenarios were computed within seconds.


\section{Results}\label{sec:results}
\begin{figure}
    \centering
    \includegraphics[width=4.5in]{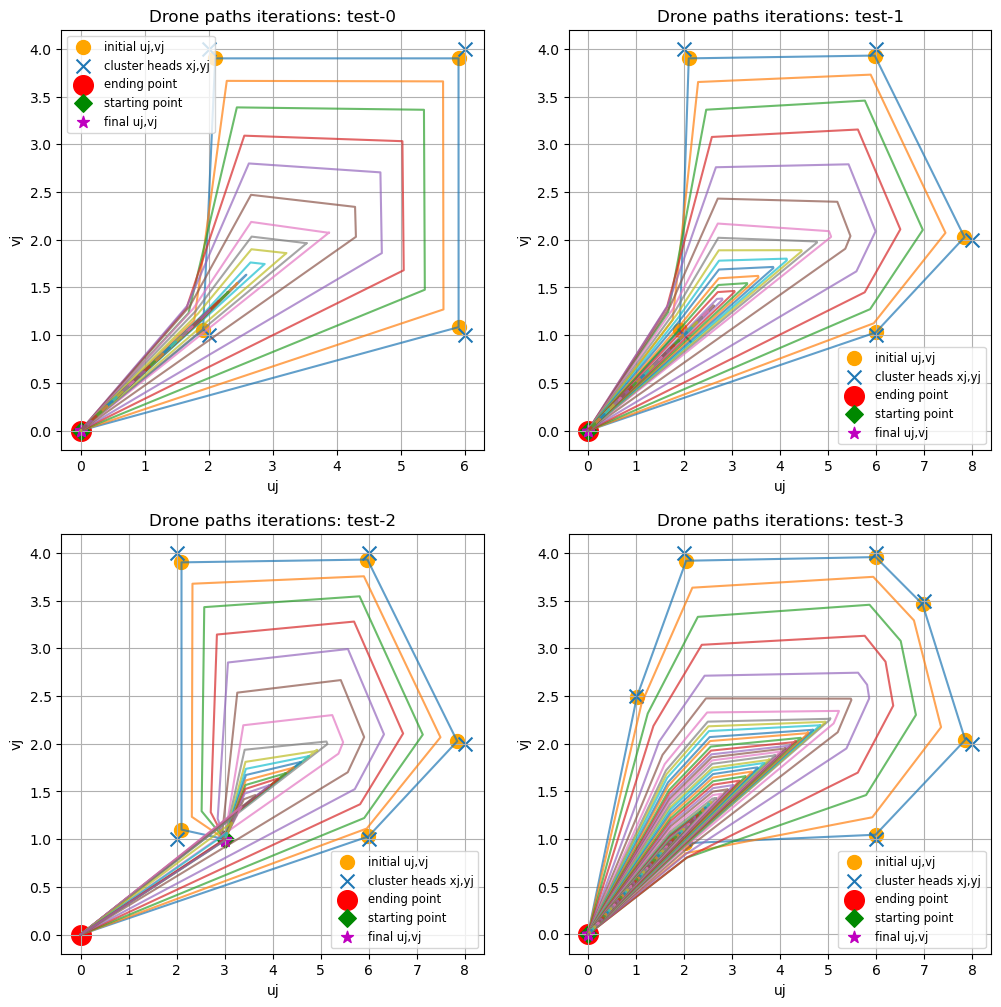}
    \caption{Drone path homotopy solutions  for 4 different test scenarios.}
    \label{fig:chart-1}
\end{figure}

\begin{figure}
    \centering
    \includegraphics[width=4.5in]{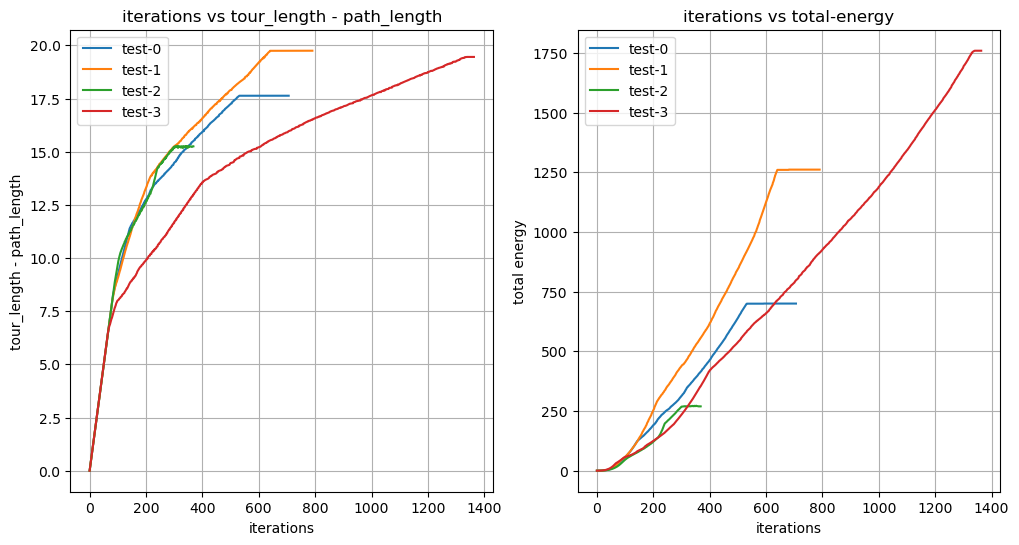}
    \caption{\emph{(left)} Iteration number versus path defect  and \emph{(right)} iteration number versus total energy for the four test scenarios shown in Figure~\ref{fig:chart-1}.}
    \label{fig:chart-2}
\end{figure}

\begin{figure}
    \centering
    \includegraphics[width=4.5in]{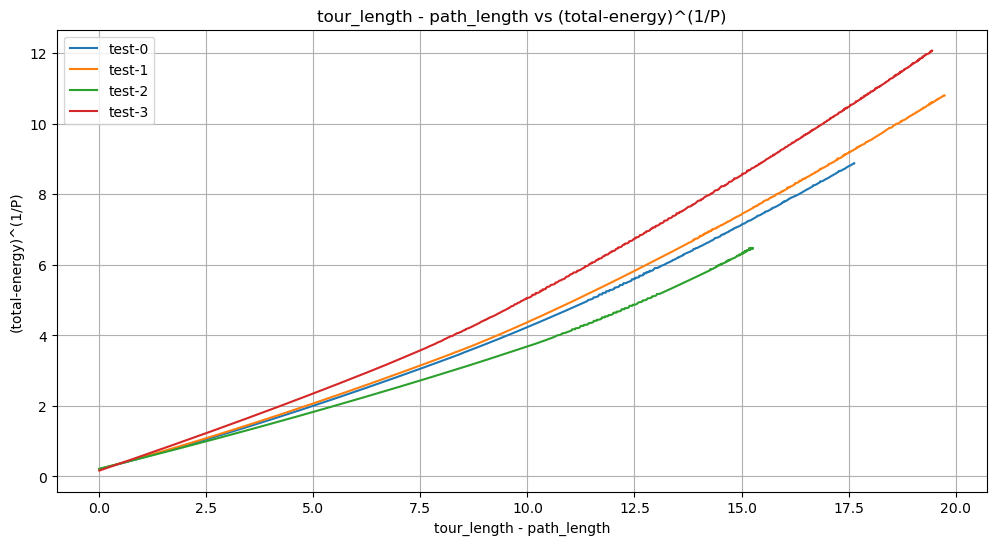}
    \caption{Path defect vs $\text{(total energy)}^{1/P}$}
    \label{fig:chart-3}
\end{figure}

 Figures~\ref{fig:chart-1}-\ref{fig:chart-3} display results from simulations of the test cases described in Section~\ref{sec:testCases} using the python code described in Section~\ref{sec:python}.  To simplify the discussion, we will define the \emph{path defect} as the difference between the length of the minimum tour joining the cluster heads and the length of the drone's path.

Figure~\ref{fig:chart-1} contains a $2\times 2$ grid of plots that shows several results from each of the above test scenarios. The chart correspond to the drone paths generated by the homotopy solution, starting from the traveling-salesman path that joins the cluster heads and decreasing down to a minimum straight-line path that joins starting and ending points. We can see that initially the paths roughly preserve the original shape as they shrink. (In fact, it can be shown that as the path shrinks the path vertices always move in the direction of the angle bisectors of the path polygon.) 
After several iterations, path vertices merge: for example, the first figure in the first row starts with a path containing 5 segments, which eventually becomes 4 and then 3 segments as the path length shrinks.  Similar changes occur in the other three figures.

The left-hand graph in Figure~\ref{fig:chart-2} shows the path defect  as a function of step number for each of the 4 test cases. The initial path length is equal to the tour length and hence the graphs start from 0. All graphs start out with the same upward slope, which is equal to 1/(step size), since $dg/ds = -1$ according to Eq.\eqref{eq:gDeriv}. 
However, when the drone's information-gathering points merge, the slope changes. This indicates that there may be numerical problems with the solution following point mergers, because eq. \eqref{eq:gDeriv} should guarantee a constant slope of 1/(step size).  
This phenomenon will be further investigated in the future (see Section~\ref{sec:conclusions}

The second chart in the middle graph gives cluster head's power transmission loss as a function of step number. As expected, this increases as the path length decreases.

The tradeoff  between drone path length and transmission loss is more clearly shown in Figure~\ref{fig:chart-3}, which plots the cluster head energy consumption on the vertical axis versus  cluster head tour length minus drone path length on the horizontal axis. In all cases, the $p$'th root of energy initially depends nearly linearly on path defect. The $p$'th root of energy is the $\ell^p$ norm of the  vector $[|\vec{z}_1-\vec{w}_1|, \ldots, |\vec{z}_J-\vec{w}_J| ]$ consisting of distances between cluster heads and corresponding drone path vertices.  In the case where all distances $|\vec{z}_j-\vec{w}_j|$ are approximately equal, then the $\ell^p$ norm is nearly proportional to the mean distance between cluster head and corresponding drone vertex. The linear section of the graphs in Figure~\ref{fig:chart-3} indicate an initial linear dependence of the distance between cluster heads and drone vertices as the path defect increases. However, as the drone range continues to decrease we see an increase in the rate of  consistent convex dependence 


\clearpage

\section{Conclusions}\label{sec:conclusions}

The solution described above for
finding optimal drone information-harvesting paths is locally optimal, and globally optimal for sufficiently long drone ranges. The numerical solution algorithm executes quickly enough that it can be implemented in real time. The solution can be used either to optimize total transmit power consumption, and maximum power consumption by cluster heads. However, there are limitations in that the speed of execution slows considerably when the drone range decreases below a certain point (i.e. when path vertices begin to merge), and the characteristics of the numerical solution indicate that the solution may no longer be optimal.

\section{Future work}\label{sec:future}

The following  future work is proposed in order to improve the solution.
\begin{itemize}
    \item  Further investigations may be made into the algorithm's behavior after drone vertex mergers. One possibility is to modify the algorithm so that when path vertices merge, the two associated cluster heads may be replaced by a single virtual cluster head that produces the same power loss.  Then the homotopy solution can be continued with vectors $\vec{u},\vec{v}$ that each have one less component.
    \item The algorithm can be modified to accommodate moving cluster heads, by making $\vec{x}$ and $\vec{y}$ in  \eqref{eq:def_f} to be functions  of $\vec{u}$ and $\vec{v}$.
    \item Further explorations of the question of global optimization may be pursued. different homotopies associated with different orderings of the cluster heads.
    \item The numerical homotopy solution code can be improved in a number of ways to decrease execution times. In particular, existing solver packages and/or variable step-size solution algorithms may be explored.
\end{itemize}

\section{Appendix: Python implementation}\label{sec:python}
\subsection{Link to code}
The python code used to generate the figures shown in Section~\ref{sec:results} is available at: \url{https://github.com/ganap-ram/drone}.
\subsection{Implemented equations}
This section describes the differential equations used to solve  as implemented in the python code. Although the implementation is mathematically equivalent to the description in Section~\ref{sec:homotopy}, the notation used is somewhat different.

For notational and coding simplicity we define some intermediate variables:
\begin{equation*}
\begin{aligned}
{u_j}, {v_j} &= j^{th} \textrm{ turning point of drone, }j = {1..J}\\
{a_j} &= {u_j - x_j}\\
{b_j} &= {v_j - y_j}\\
A_j &= {a_j}^2 + {b_j}^2 \\
{m_j} &= u_{j} - u_{j+1}\\
{n_j} &= v_{j} - v_{j+1}\\
M_j &= {m_j}^2 + {n_j}^2 \\
H_j &= {M_j}^{3/2} \\
q &= \frac{p}{2} - 1\\
r &= 1/2\\
\end{aligned}    
\end{equation*}

The following are the various coefficient values of the derivatives in the matrix:
\begin{equation}
\begin{aligned}    
s_{j1} &= \lambda\frac{M_{j-1} - m_{j-1}^{2}}{H_{j-1}}\\
s_{j2} &= -\lambda\frac{m_{j-1}n_{j-1}}{H_{j-1}}\\
s_{j3} &= 2pqa_{j}^{2}A_{j}^{q-1} + pA_{j}^{q} - \lambda\frac{H_{j}(M_{j-1} - m_{j-1}^{2}) + H_{j-1}(M_{j} - m_{j}^{2})}{H_{j-1}H_{j}}\\
s_{j4} &= 2pqa_{j}b_{j}A_{j}^{q-1} + \lambda\frac{m_{j-1}n_{j-1}H_{j} + m_{j}n_{j}H_{j-1}}{H_{j-1}H_{j}}\\
s_{j5} &= \lambda\frac{M_{j} - m_{j}^{2}}{H_{j}}\\
s_{j6} &= -\lambda\frac{m_{j}n_{j}}{H_{j}}\\
w_{j} &= \frac{-m_{j-1}\sqrt{M_{j}} + m_{j}\sqrt{M_{j-1}}}{\sqrt{M_{j-1}M_{j}}}
\end{aligned}
\end{equation}

\begin{equation}
    \begin{aligned}
t_{j1} &= -\lambda\frac{m_{j-1}n_{j-1}}{H_{j-1}}\\
t_{j2} &= \lambda\frac{M_{j-1} - n_{j-1}^{2}}{H_{j-1}}\\
t_{j3} &= 2pqa_{j}b_{j}A_{j}^{q-1} + \lambda\frac{m_{j-1}n_{j-1}H_{j} + m_{j}n_{j}H_{j-1}}{H_{j-1}H_{j}}\\
t_{j4} &= 2pqb_{j}^{2}A_{j}^{q-1} + pA_{j}^{q} - \lambda\frac{H_{j}(M_{j-1} - n_{j-1}^{2}) + H_{j-1}(M_{j} - n_{j}^{2})}{H_{j-1}H_{j}}\\
t_{j5} &= -\lambda\frac{m_{j}n_{j}}{H_{j}}\\
t_{j6} &= \lambda\frac{M_{j} - n_{j}^{2}}{H_{j}}\\
z_{j} &= \frac{-n_{j-1}\sqrt{M_{j}} + n_{j}\sqrt{M_{j-1}}}{\sqrt{M_{j-1}M_{j}}}
\end{aligned}
\end{equation}

And the differential equations have the form for $j=1\ldots J$ (note that $\frac{du_{0}}{ds}=\frac{dv_{0}}{ds}=0$, so that the first two equations below have 4 terms instead of 6 when $j=1$):
\begin{equation}\label{eq:DEs}
\begin{aligned}
& s_{j1}\frac{du_{j-1}}{ds} + s_{j2}\frac{dv_{j-1}}{ds} + s_{j3}\frac{du_{j}}{ds} + s_{j4}\frac{dv_{j}}{ds} + s_{j5}\frac{du_{j+1}}{ds} + s_{j6}\frac{dv_{j+1}}{ds} - w_j\frac{d\lambda}{ds} &= 0\\
& t_{j1}\frac{du_{j-1}}{ds} + t_{j2}\frac{dv_{j-1}}{ds} + t_{j3}\frac{du_{j}}{ds} + t_{j4}\frac{dv_{j}}{ds} + t_{j5}\frac{du_{j+1}}{ds} + t_{j6}\frac{dv_{j+1}}{ds} - z_j\frac{d\lambda}{ds} &= 0\\
\frac{dg}{ds} &= -1 \implies \sum_{j=1}^J \left(\frac{\partial g}{\partial u_j}\frac{du_j}{ds} + \frac{\partial g}{\partial u_j}\frac{dv_j}{ds}\right) &= -1,
\end{aligned}
\end{equation}

This is represented in the matrix notation as:
\begin{equation*}
\begin{aligned}
K D  = C,
\end{aligned}
\end{equation*}
where:
\begin{align}
K &=  
\begin{bmatrix}
\begin{matrix}
s_{13} & s_{14} & s_{15} & s_{16} & 0 & 0 & 0 & 0 & \cdots & 0 & 0 & -w_{1}\\ 
t_{13} & t_{14} & t_{15} & t_{16} & 0 & 0 & 0 & 0 & \cdots & \cdots & 0 & -z_{1}\\ 
s_{21} & s_{22} & s_{23} & s_{24} & s_{25} & s_{16} & 0 & 0 & \cdots & \cdots & 0 & -w_{2}\\
t_{21} & t_{22} & t_{23} & t_{24} & t_{25} & t_{16} & 0 & 0 & \cdots & \cdots & 0 & -z_{2}\\
0 & 0 & s_{31} & s_{32} & s_{33} & s_{34} & s_{35} & s_{36} & \cdots & \cdots & \cdots & \vdots\\
0 & 0 & t_{31} & t_{32} & t_{33} & t_{34} & t_{35} & t_{36} & \cdots & \cdots & \cdots & \vdots\\
0 & 0 & \cdots & \cdots & \cdots & \cdots & \cdots & \cdots & \cdots & \cdots &\vdots & \vdots\\
0 & 0 & \cdots & \cdots & \cdots & \cdots & \cdots & \cdots & \cdots & \cdots &0  & \vdots\\
0 & 0 &\cdots & \cdots & \cdots & 0 & 0 & s_{J1} & s_{J2} & s_{J3} & s_{J4} & -w_{J}\\
0 & 0 & 0 & \cdots & \cdots & 0 & 0 & t_{J1} & t_{J2} & t_{J3} & t_{J4}& -z_{J}\\
w_{0} & z_{0} & w_{1} & z_{1} & \cdots & \cdots & \cdots & \cdots & \cdots & w_{J} & z_{J}& 0\\
\end{matrix}
\end{bmatrix}\\
\end{align}
\begin{align}
D &= 
\begin{bmatrix}
{du}_{1}/{ds} \\
{dv}_{1}/{ds} \\
{du}_{2}/{ds} \\
{dv}_{2}/{ds} \\
{du}_{3}/{ds} \\
{dv}_{3}/{ds} \\
\vdots \\
\vdots \\
\vdots \\
{du}_{J}/{ds} \\
{dv}_{J}/{ds} \\
{d\lambda}/{ds} \\
\end{bmatrix}\\
\end{align}
And
\begin{align}
C &= 
\begin{bmatrix}
0\\
0\\
0\\
0\\
0\\
0\\
\vdots \\
\vdots \\
\vdots \\
0\\
0\\
-1\\
\end{bmatrix}
\end{align}
\section{Acknowledgements}
Thanks to Dr. Aristide Tsemo for discovering several errors in the equations in Section~\ref{sec:derivation}.

\printbibliography

\end{document}